\documentclass[12pt]{article}
\input{epsf}

\usepackage[english]{babel}
\usepackage{amsmath}
\usepackage{amsthm,amssymb,latexsym}
\usepackage{graphics,psfig}
\newcommand{\psdraw}[3]{\begin{array}{c} \hspace{-1mm}
\raisebox{-4pt}{\psfig{figure=#1.ps,width=#2,height=#3}}
\hspace{-245pt}\end{array}}

\newtheorem{The}{Theorem}

\newtheorem{Exa}{Example}
\newtheorem{Rem}{Remark}

\newtheorem{Obs}{Observation}
\newtheorem{Pro}{Proposition}

\def\aa{$\psdraw{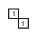}{90mm}{7mm}$}

\def\bb{$\psdraw{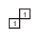}{90mm}{7mm}$}

\def\h1{$\psdraw{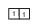}{85mm}{4mm}$\hspace{15pt}}

\def\v1{$\psdraw{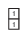}{85mm}{7mm}$\hspace{8pt}}
\begin{document}
\title {\sc the problem of the pawns}
\author {Sergey Kitaev and Toufik Mansour \footnote{Research
financed by EC's IHRP Programme, within the Research Training
Network "Algebraic Combinatorics in Europe", grant
HPRN-CT-2001-00272}} \maketitle
\begin{center}
Matematik, Chalmers tekniska h\"ogskola och G\"oteborgs
universitet, 412~96  G\"oteborg, Sweden\\

{\it kitaev@math.chalmers.se, toufik@math.chalmers.se}
\end{center}

\def\mn{\mbox{-}}
\def\SS{{\mathcal S}}
\newcommand\mps[1]{\marginpar{\small\sf#1}}
\section*{Abstract}
In this paper we study the number $M_{m,n}$ of ways to place
nonattacking pawns on an $m\times n$ chessboard. We find an upper
bound for $M_{m,n}$ and analyse its asymptotic behavior. It turns
out that $\lim_{m,n\rightarrow\infty}(M_{m,n})^{1/mn}$ exists and is
bounded from above by $(1+\sqrt{5})/2$. Also, we consider a lower bound for
$M_{m,n}$ by reducing this problem to that of tiling an $(m+1)\times (n+1)$
board with square tiles of size $1\times 1$ and $2\times 2$. Moreover,
we use the transfer-matrix method to implement an algorithm that
allows us to get an explicit formula for $M_{m,n}$ for given $m$.

\medskip
\noindent {\sc 2000 Mathematics Subject Classification}: 05A16,
05C50, 52C20, 82B20
\section{Introduction}
On an $m\times n$ chessboard, we place a number of nonattacking
pawns, all of the same colour, say white. The main question here
is: How many different placements are possible? A similar problem
concerning placements of the maximum number of nonattacking kings
on a $2m\times 2n$ chessboard is treated in~\cite{Wilf}. The main
result of that paper is the following theorem.

\begin{The}
Let $f_m(n)$ denote the number of ways that $mn$ nonattacking
kings can be placed on a $2m\times 2n$ chessboard. For each
$m=1,2,3,\ldots $ there are constants $c_m>0$, $d_m$, and $0\le
\theta_m<m+1$ such that
$$f_m(n)= (c_mn+d_m)(m+1)^n+O(\theta_m^n)\qquad (n\to\infty).$$
\end{The}

Given an $m\times n$ chessboard. We mark a square containing a pawn
by 1, and a square that does not contain a pawn by~0. The placement
of pawns is then completely specified by an $m\times n$ binary
matrix. Moreover, to be a legal placement, the binary matrix
cannot contain the following two letter words: \aa and \bb (here
we use the fact that all pawns are of the same colour and thus
they are allowed to attack at the same directions: either at
North-West and North-East or at South-West and South-East). For
example, the matrix
$$
\begin{array}{cccccc}
1 & 0 & 1 & 1 & 0 & 1\\
1 & 0 & 0 & 0 & 0 & 0\\
0 & 0 & 1 & 0 & 1 & 1\\
\end{array}
$$
corresponds to a legal placement of pawns on a $3\times6$ board.
So, our main question can be reformulated as follows: How many
binary $m\times n$ matrices simultaneously avoid the words \aa and
\bb? We denote the number of such matrices by $M_{m,n}$.

Studying matrices avoiding certain words, and thus studying our
original problem, is interesting, for instance, from a graph
theoretic point of view~\cite{CalkinWilf}. In that paper, the
authors considered the (vertex) independence number of the
$m\times n$ grid graph using the matrices with the property that
no two consecutive 1's occur in a row or a column.

In this paper, we use the transfer-matrix approach to implement
an algorithm that allows us to find a formula for $M_{m,n}$ for any
given $m$ (see Sections~\ref{tilings} and~\ref{transfer}).
Moreover, in Section~\ref{upper} we find an upper bound for $M_{m,n}$
and, in Section~\ref{tilings}, we discuss
how the tiling problem is related to finding a lower bound for $M_{m,n}$.
Also, we prove that the double limit
$\lim_{m,n\rightarrow\infty} (M_{m,n})^{1/mn}$ exists and is bounded
from above by $(1+\sqrt{5})/2$ (see Sections~\ref{transfer}
and~\ref{upper}). Finally, in Section~\ref{formulas}, we suggest an
approach to study $M_{m,n}$, which, in particular, allows to prove that
$M_{2m,n}$ is a perfect square (see Theorem~\ref{perfectSquare}). Using
this approach we obtain formulas for $M_{m,n}$, where $2\leq m\leq 6$.

\section{The upper bound for $M_{m,n}$}
\label{upper}
To obtain an upper bound for $M_{m,n}$, we determine the
number $U_{m,n}$ of binary $m\times n$ matrices that avoid the
word \aa. Of course, $U_{m,n}$ counts also the number of binary
$m\times n$ matrices that avoid the word \bb, which follows from
arranging the columns of all matrices under consideration in reverse order
(in particular, \aa is the reverse of \bb).

The following theorem gives a formula for the number of binary
matrices that avoid the word \aa in terms of the Fibonacci
numbers.

\begin{The}\label{about_U_m_n}
For any $n,m\geq0$,
$$U_{m,n}=\left\{
\begin{array}{ll}
F_{m+1}^{n-m+1}
\left(\displaystyle\prod_{i=0}^{m}F_i\right)^2,& \mbox{ if } n\geq m,\\[3mm]
F_{n+1}^{m-n+1}
\left(\displaystyle\prod_{i=0}^{n}F_i\right)^2,& \mbox{ if } n<m,
\end{array}\right.
$$
where $F_i$ is the $i$-th Fibonacci number defined by $F_0=F_1=1$, and
$F_{n+2}=F_{n+1}+F_n$ for $n\geq 0$.
\end{The}

\begin{proof}
Let $A$ be an $m\times n$ (0,1)-matrix that avoids the word \aa. We
change the shape of $A$ using the following procedure. We shift
the first column of $A$ one position down with respect to the second
column. In the obtained shape, we shift the first and second
columns one position down with respect to the third column, and so
on. After shifting with respect to the $n$-th column, one obtains
the shape $\bar{A}$, that has the form similar to that on Figure~\ref{fig1}.

\begin{center}
\begin{figure}[h]
$$\hspace{130pt}\psdraw{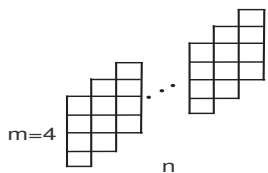}{180mm}{23mm}\hspace{-14pt}$$
\caption{The shape $\bar{A}$ for $m=4$.} \label{fig1}
\end{figure}
\end{center}

It is easy to see that $A$ avoids the word \aa if and only if
$\bar{A}$ avoids the word \h1. But $\bar{A}$ avoids \h1 if and
only if each row of $\bar{A}$ avoids \h1 (there are no additional
restrictions). This is well known and is not difficult to see that
the number of different binary strings of length $\ell$ that avoid
\h1 is given by $F_{\ell+2}$.

To find $U_{m,n}$, it remains to find out the lengths of
the rows in $\bar{A}$, and since these rows are independent from
each other, to multiply together the corresponding Fibonacci
numbers. If $n\geq m$, $\bar{A}$ has two rows of each of the
following lengths: $1$, $2, \ldots, m-1$, and $n-m+1$ rows of
length $m$. So, in this case
$$U_{m,n}=F_{m+1}^{n-m+1}
\left(\displaystyle\prod_{i=0}^{m}F_i\right)^2.$$
The case $n<m$
is given by changing $m$ by $n$, and $n$ by $m$ in the
considerations above.
\end{proof}

\begin{center}
\begin{figure}[h]
$$\hspace{10pt}\psdraw{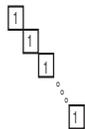}{110mm}{18mm}$$
\caption{The $k$-diagonal word.} \label{fig2}
\end{figure}
\end{center}

Let $A$ be any binary matrix, we say that $A$ avoids the {\em
$k$-diagonal word} (see Figure~\ref{fig2}) if there are no $k$
consecutive 1's in any diagonal of $A$. Theorem~\ref{about_U_m_n}
can be generalized to the case of avoiding the $k$-diagonal word.
This generalization involves the $k$-generalized Fibonacci
numbers. We do not use the generalization to proceed with the
problem of the pawns, but we state it as Theorem~\ref{gener_U_m_n}
because it is interesting by its own.

Let $F_{k,n}$ be the {\em $n$-th $k$-generalized Fibonacci
number} defined by $F_{k,n}=0$ for $n<0$, $F_{k,0}=1$, and
$F_{k,n}=F_{k,n-1}+F_{k,n-2}+\cdots +F_{k,n-k}$ for $n\geq1$ (for
example, see~\cite{F,SP}).

Let $U_{m,n}(k)$ denote the number of $m\times n$ binary matrices
that avoid the $k$-diagonal word. The following theorem can be
proved using the same arguments as those in
Theorem~\ref{about_U_m_n} and the observation that the number of
different binary strings of length $\ell$ that avoid the word
$\underbrace{11\ldots 1}_{k}$ is given by $F_{k,\ell+1}$ (we leave
this observation as an exercise).

\begin{The}\label{gener_U_m_n}
Let $k\geq 2$. For all $n,m\geq0$,
$$U_{m,n}(k)=\left\{
\begin{array}{ll}
F_{k,m+1}^{n-m+1}
\left(\displaystyle\prod_{i=0}^{m}F_{k,i}\right)^2,& \mbox{ if } n\geq m,\\[3mm]
F_{k,n+1}^{m-n+1}
\left(\displaystyle\prod_{i=0}^{n}F_{k,i}\right)^2,& \mbox{ if }
n<m,
\end{array}\right.
$$
where $F_{k,i}$ is the $i$-th $k$-generalized Fibonacci number.
\end{The}

As a corollary to Theorem~\ref{about_U_m_n}, we get an upper bound for
$M_{m,n}$. Indeed, $M_{m,n}\leq U_{m,n}$ since $U_{m,n}$ deals with
avoidance of \aa, whereas $M_{m,n}$ deals additionaly with one more
restriction, namely \bb. We state this result as the following theorem.

\begin{The}\label{upper_bound}
We have
$$M_{m,n}\leq
\left\{
\begin{array}{ll}
F_{m+1}^{n-m+1}
\left(\displaystyle\prod_{i=0}^{m}F_i\right)^2,& \mbox{ if } n\geq m,\\[3mm]
F_{n+1}^{m-n+1}
\left(\displaystyle\prod_{i=0}^{n}F_i\right)^2,& \mbox{ if } n<m,
\end{array}\right.$$
where $F_i$ is the $i$-th Fibonacci number.
\end{The}

The upper bound for $M_{m,n}$ involves the product of the first nonzero
Fibonacci numbers. It is known~\cite[A003266]{SP}
that an asymptotic for the product of the first $n$ nonzero Fibonacci
numbers is given by
\begin{equation}\label{productFib}
\frac{c}{\sqrt{5}^{n-1}}\left(\frac{1+\sqrt{5}}{2}\right)^{\frac{n(n-1)}{2}},
\end{equation}
where
$c=\prod_{j\geq1}\left(1-\left(\frac{\sqrt{5}-3}{2}\right)^j\right)=1.2267420107203532444176302\cdots$.
This result and Theorem~\ref{upper_bound} give the following
theorem.

\begin{The}
We have
$$\lim_{n,m\rightarrow\infty}(M_{m,n})^{\frac{1}{mn}}\leq\frac{1+\sqrt{5}}{2}.$$
\end{The}
\begin{proof}
The existence of the limit
$\lim_{n,m\rightarrow\infty}(M_{m,n})^{\frac{1}{mn}}$ is proved in
Theorem~\ref{thexists}. Using Theorem~\ref{upper_bound}, it is
enough to prove that
\begin{equation}\label{fN}
\lim_{n,m\rightarrow\infty}(U_{m,n})^{\frac{1}{mn}}=\frac{1+\sqrt{5}}{2}.
\end{equation}
For given two functions $f(n)$ and $g(n)$, we define $f(n)\sim g(n)$ if
$\lim\limits_{n\rightarrow\infty}\frac{f(n)}{g(n)}=1$.

Suppose $n\geq m$. By~(\ref{productFib}) we have
$$\left(\prod_{i=0}^{m}F_i\right)^2\sim \frac{c^2}{\sqrt{5}^{2m-2}}
\left(\frac{1+\sqrt{5}}{2}\right)^{m^2-m},$$
and using the formula for the Fibonacci numbers, namely
$$F_m=\frac{1}{\sqrt{5}}\left(\left(\frac{1+\sqrt{5}}{2}\right)^{m+1}-\left(\frac{1-\sqrt{5}}{2}\right)^{m+1}\right),$$ we obtain that
$$U_{m,n}\sim \frac{c^2}{\sqrt{5}^{n+m-1}}
\left(\frac{1+\sqrt{5}}{2}\right)^{nm+2n-2m+2}.$$
This formula holds for the case $n<m$, by replacing $m, n$ by $n, m$ in
the considerations above. Hence, (\ref{fN}) holds.
\end{proof}

\section{Tiling rectangles and a lower bound for $M_{m,n}$}\label{tilings}

Let $L_{m,n}$ denote the number of $m\times n$ binary matrices that
simultaneously avoid the words \aa, \bb, \v1 and \h1. Clearly,
$L_{m,n}\leq M_{m,n}$, since when we deal with $L_{m,n}$ we have more
restrictions than when we consider $M_{m,n}$.
Thus, we are interested in finding the numbers $L_{m,n}$,
that give us a lower bound for $M_{m,n}$. In this section we show
that $L_{m,n}$, in fact, gives the number of tilings of an $(m+1)\times (n+1)$
area with square tiles of size $1\times 1$ and $2\times 2$ which was
studied in~\cite{Heubach} and~\cite{ChiHeub}. So, the number of the
tilings is equal to the number of $m\times n$ binary matrices that
avoid the words \aa, \bb, \v1 and \h1. A bijection $\theta$ between these two
combinatorial objects is given by the following.

Let $A$ be an $m\times n$ matrix that avoids the words \aa, \bb, \v1 and
\h1. We make from $A$ an $(m+1)\times (n+1)$ matrix $\bar{A}$ by adjoin an
additional
$m\times 1$ column consisting of 0's from the right side, and an additional
$1\times (n+1)$ row, also having only 0's, from below. Now, once we
meet an occurrence of $1$ in $\bar{A}$, we place a $2\times 2$ tile in
such way, that the $1$ appears in the top-left corner
of the tile. After considering all 1's and placing corresponding $2\times 2$
tiles, we fill in the uncovered squares of $\bar{A}$ by $1\times 1$ tiles. The
fact that $A$ avoids \aa, \bb, \v1 and \h1 guarantees that covering in the
way proposed by us is non-overlapping, and thus we get a tiling of
an $(m+1)\times (n+1)$ board.

Conversely, for any given tiling with square tiles of size $1\times 1$ and
$2\times 2$, we can place 1 in the top-left corner of any $2\times2$ tile,
0's in the other squares, and remove the rightmost column and the bottom row.
Obviously, we get an $m\times n$ binary matrix that avoids the words
\aa, \bb, \v1 and \h1.

\begin{center}
\begin{figure}[h]
\hspace*{20pt}
\epsfxsize=350.0pt
\epsffile{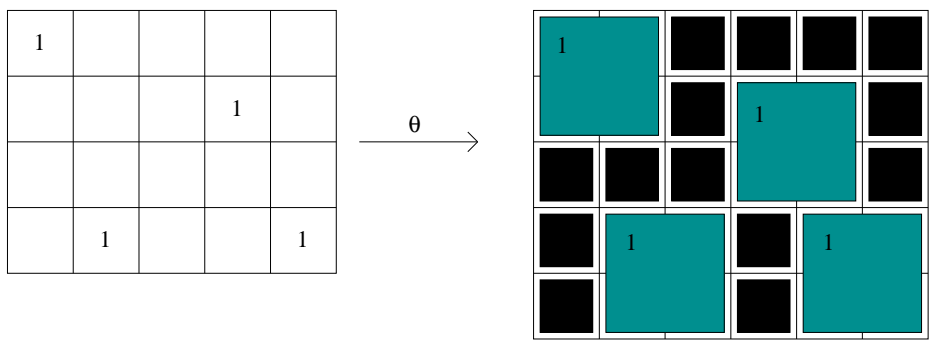}
\caption{The bijection $\theta$.} \label{fig3}
\end{figure}
\end{center}

Figure~\ref{fig3} shows how the bijection $\theta$ works in the case of
a $4\times5$ matrix.

Unfortunately, we cannot get much use of the papers~\cite{Heubach}
and~\cite{ChiHeub}, since there, for our purpose, one has explicit
formulas only for $m=2,3$, and for $m=4,5$ one has recursive
formulas only. That means that from that source, we have an information about
$L_{m,n}$, where $1\leq m\leq 4$.
\begin{Rem}\label{lf}
If we use the transfer-matrix approach {\rm(see~\cite[pawns-kings]{M})}
for finding the formula for $L_{m,n}$, where $m$
is given, then we get the following: $L_{1,n}=F_{n+1}$,
$L_{2,n}=\frac{1}{3}(2^{n+2}-(-1)^n)$, and
$$\begin{array}{l}
L_{3,n}=a\left(\frac{2}{3}+\frac{2\sqrt{13}}{3}\cos\beta\right)^n
+b\left(\frac{2}{3}-\frac{2\sqrt{39}}{3}\sin\beta-\frac{2\sqrt{13}}{3}\cos\beta\right)^n\\
\qquad\qquad\qquad\qquad\qquad\qquad\qquad\qquad
+c\left(\frac{2}{3}+\frac{2\sqrt{39}}{3}\sin\beta-\frac{2\sqrt{13}}{3}\cos\beta\right)^n,
\end{array}$$
where
$\beta=\frac{1}{3}\arctan\left(\frac{3}{8}\sqrt{237}\right)$,
$a\approx1.51212496094$, $b\approx-0.542960193686$, and
$c\approx0.0308352327442$.
\end{Rem}

\section{The transfer-matrix method}\label{transfer}
We use the transfer-matrix method, in a manner that is similar to way it was
used in~\cite{Wilf}.

For $m$, $n$ fixed, we can think of constructing the $m\times n$ binary
matrices avoiding the words \aa and \bb by gluing together columns
that are chosen from the collection of possible columns, making sure that
when we glue an additional column onto the right-hand edge of the
structure, the new column does not come into conflict with the previous
right-hand column. The collection of possible columns $C_m$ is the set
of all $m$-vectors $v$ of 0's and 1's. Clearly, $|C_m|=2^m$.

The condition that vectors $v$, $w$ in $C_m$ are possible consecutive
pair of columns in a matrix avoiding \aa and \bb is simply that
$v_iw_{i+1}=0$ and $v_{i+1}w_i=0$ for all $i$, $1\leq i\leq m-1$. We say
that such $v$ and $w$ are {\em cross-orthogonal}.

Thus, all possible matrices under consideration are obtained by beginning
with some vector
of $C_m$, and in general, having arrived at some sequence of vectors of
$C_m$, adjoin any vector of $C_m$ that is cross-orthogonal to the last one
previously chosen until $n$ vectors have been selected.

We define a matrix $T=T_m$, the transfer matrix of the problem, as follows.
$T$ is an $2^m\times 2^m$ symmetric matrix of $0$'s and $1$'s whose rows
and columns are indexed by vectors of $C_m$. The entry of $T$ in position
$(v,w)$ is 1 if the vectors $v$, $w$ are cross-orthogonal, and $0$
otherwise. $T$ depends only on $m$, not on $n$.

Let $M_{m,n}(u)$ denote the number of $m\times n$ binary
matrices avoiding the words \aa and \bb whose rightmost column vector is $u$.
Then, clearly, we have
$$M_{m,n+1}(v)=\sum_{u\in C_m}M_{m,n}(u)T_{u,v}\ \ \ \ (n\geq0; v\in C_m),$$
or, in matrix-vector notation, $M_{n+1}=TM_n$, with $M_0={\bf 1}$ the
vector of length $2^m$ whose entries are all 1's. It follows that
$M_n=T^n\cdot{\bf 1}$, for
all $n\geq 0$. The number of matrices $M_{m,n}$ is the sum of the
entries of the vector $M_n$. Thus, if ${\bf 1^{\prime}}$ denote the row
of length $2^m$ whose entries are all 1's, we have
$$M_{m,n}={\bf 1}^{\prime}\cdot T^n\cdot{\bf 1},$$
i.e., $M_{m,n}$ is the sum of all of the entries of the matrix $T^n$.

\begin{Exa} The transfer-matrices $T_2$ and $T_3$ {\rm(see~\cite[pawns]{M})}
are given, for instance, by
$$T_2=\begin{pmatrix}
1& 1& 1& 1\\
1& 1& 0& 0\\
1& 0& 1& 0\\
1& 0& 0& 0
\end{pmatrix}\quad\mbox{and}\quad
T_3=\begin{pmatrix}
1& 1& 1& 1& 1& 1& 1& 1\\
1& 1& 0& 0& 1& 1& 0& 0\\
1& 0& 1& 0& 0& 0& 0& 0\\
1& 0& 0& 0& 0& 0& 0& 0\\
1& 1& 0& 0& 1& 1& 0& 0\\
1& 1& 0& 0& 1& 1& 0& 0\\
1& 0& 0& 0& 0& 0& 0& 0\\
1& 0& 0& 0& 0& 0& 0& 0
\end{pmatrix}.$$
\end{Exa}

Since $T$ has nonnegative entries, its dominant eigenvector cannot
be orthogonal to ${\bf 1}$, and so we have at once that
$\lim_{n\rightarrow \infty}(M_{m,n})^{\frac{1}{n}}$ exists for
each~$m$, and is equal to $\alpha_m$, the largest eigenvalue of
the transfer-matrix $T$ (real and symmetric matrix). It follows
that
\begin{equation}\label{infsup}
\liminf_{m}(\alpha_{m})^{\frac1m}=\liminf_{m,n}(M_{m,n})^{\frac{1}{mn}}\leq
\limsup_{m,n}
(M_{m,n})^{\frac{1}{mn}}=\limsup_m(\alpha_{m})^{\frac1m}.
\end{equation}

\begin{The}\label{thexists}
The limit $\lim_{m,n\rightarrow\infty}(M_{m,n})^{\frac{1}{mn}}$
exists.
\end{The}
\begin{proof}
By the fact that $T$ is symmetric and real matrix together with
using the maximum principle we get, for any $q\geq1$,
$$\frac{({\bf1},(T_m)^q\cdot{\bf1})}{({\bf1},{\bf1})}\leq
(\alpha_m)^q.$$ Since $M_{m,q}=M_{q,m}$ by the definitions, we
have $({\bf1},(T_m)^q\cdot{\bf1})=({\bf1},(T_q)^m\cdot{\bf1})$.
Hence,
$$\left(\frac{({\bf1},(T_q)^m\cdot{\bf1})}{({\bf1},{\bf1})}\right)^{\frac1m}\leq
(\alpha_m)^\frac{q}{m}.$$ Taking the $\liminf_{m}$ of both sides
of the inequality above, together with using the fact that
$|C_m|=2^m$, we have $\frac{\alpha_q}{2}\leq
\left(\liminf_m(\alpha_m)^\frac{1}{m}\right)^q$, which implies
$$\limsup_q(\alpha_q)^{\frac{1}{q}}\leq
\liminf_m(\alpha_m)^\frac{1}{m}.$$ Using~(\ref{infsup}) we get the
desired result.
\end{proof}

Using the transfer-matrix approach one can obtain an explicit
formula for $M_{m,n}$, where $m\geq1$ is given. We
implemented an algorithm for finding the transfer-matrix $T_m$ in
Maple (see~\cite[pawns]{M}). This algorithm yields
an explicit formula for $M_{m,n}$, where $1\leq m \leq 3$ (see
Table~1). Moreover, it finds the maximum
eigenvalue of $T_m$ for given~$m$.

\begin{table}[h]
\begin{center}
\begin{tabular}{|l|l|}\hline
  $m$ & $M_{m,n}$ \\ \hline
  $1$ & $2^n$ \\ \hline\hline
&\\[-8pt]
  $2$ &
  $\frac{7}{10}(\eta_1^{2n}+\eta_2^{2n})+\frac{3\sqrt{5}}{10}(\eta_1^{2n}-\eta_2^{2n})-\frac{2}{5}(-1)^{n}$,\\[2pt]
      & \mbox{where }$\eta_1=\frac{1}{2}(1+\sqrt{5})$\mbox{ and
      }$\eta_2=\frac{1}{2}(1-\sqrt{5})$\\[2pt]\hline\hline
&\\[-8pt]
  $3$ &
  $\frac{1}{13}\bigl(\eta_1^{n+2}+\eta_2^{n+2}\bigr)+\frac{\sqrt{3}^{n+1}}{13}\bigl(4-\sqrt{3}-(4+\sqrt{3})(-1)^n\bigr)$,\\[2pt]
      &\mbox{where }$\eta_1=\frac{1}{2}(5+\sqrt{13})$\mbox{ and
      }$\eta_2=\frac{1}{2}(5-\sqrt{13})$\\[2pt]\hline
\end{tabular}
\caption{Explicit formula for $M_{m,n}$ where $m=1,2,3$.}
\end{center}
\end{table}
For example, the maximum
eigenvalue of $T_m$ is $2$, $\left(\frac{1+\sqrt{5}}{2}\right)^2$,
$\frac{5+\sqrt{13}}{2}$ and
$\frac{8}{3}+\frac{4}{3}\sqrt{7}\cos\left(\frac{1}{3}\arctan\left(\frac{3}{67}\sqrt{111}\right)\right)$,
for $m=1,2,3,4$; respectively.

\begin{Rem}
In the case of $m=4$, the eigenvalues of $T$ are given by
$$\begin{array}{l}
\lambda_1=\frac23-\frac{4}{3}\cos\left(\frac{1}{3}\pi-\beta\right)-\frac43\sqrt{3}\sin\left(\frac13\pi-\beta\right),\\[2pt]
\lambda_2=\frac23-\frac{4}{3}\cos\left(\frac{1}{3}\pi-\beta\right)+\frac43\sqrt{3}\sin\left(\frac13\pi-\beta\right),\\[2pt]
\lambda_3=\frac83-\frac23\sqrt{7}\cos\gamma-\frac23\sqrt{21}\sin\gamma,\\[2pt]
\lambda_4=\frac83-\frac23\sqrt{7}\cos\gamma+\frac23\sqrt{21}\sin\gamma,\\[2pt]
\lambda_5=-\frac23-\frac43\cos\beta-\frac43\sqrt{3}\sin\beta,\\[2pt]
\lambda_6=-\frac23-\frac43\cos\beta+\frac43\sqrt{3}\sin\beta,\\[2pt]
\lambda_7=\frac23+\frac83\cos\left(\frac13\pi-\beta\right),\\[2pt]
\lambda_8=-\frac23+\frac83\cos\left(\frac13\pi-\beta\right),\\[2pt]
\lambda_9=\frac83+\frac43\sqrt{7}\cos\gamma,
\end{array}$$
where $\beta=\frac13\arctan\left(\frac35\sqrt{111}\right)$ and
$\gamma=\frac13\arctan\left(\frac3{67}\sqrt(111)\right)$.
\end{Rem}

\section{Formulas for $M_{m,n}$}\label{formulas}
In this section we suggest another approach to study $M_{m,n}$. In particular,
we obtain formulas for $M_{m,n}$, where $2\leq m\leq 6$ (the cases $m=2,3$
already appear in Table~1). We show how to use the following simple observation in order to investigate $M_{m,n}$.

\begin{Obs}\label{observ}
A pawn placed on a square of a chessboard cannot attack a square of
the different colour.
\end{Obs}

According to Observation~\ref{observ}, $M_{m,n}=B_{m,n}\cdot W_{m,n}$, where
$B_{m,n}$ (resp. $W_{m,n}$) is the number of ways to place nonattacking
pawns on the black (resp. white) squares of an $m\times n$ chessboard. Thus,
the original problem of finding $M_{m,n}$ can be reduced to considering
independently two shapes: that consisting of all the black squares, and
the shape consisting of all the white squares. We use this idea in the
proofs of the following theorem and propositions.

\begin{The}\label{perfectSquare} We have that $M_{2m,n}=a^2$ for some natural number $a$,
that is $M_{2m,n}$ is a perfect square.
\end{The}

\begin{proof} Using the discussion right above this theorem, it is
enough to prove that $B_{2m,n}=W_{2m,n}$. Indeed, on a $2m\times n$
chessboard, the number of black squares is the same as that of white
squares. Moreover, if we consider the shape that, say, the white squares form,
reverse it horizontally (that is, draw the rows in reverse order),
then we get exactly the same shape that the black
squares form. Also, it is easy to see that a placement of pawns before
the reversion
is legal if and only if it is legal after the reversion. Thus, we have
$B_{2m,n}=W_{2m,n}$.
\end{proof}

\begin{Pro} We have
$$M_{2,n}=(F_{n+2})^2,$$
where $F_n$ is the $n$-th Fibonacci number defined by $F_0=F_1=1$, and
$F_{n+2}=F_{n+1}+F_n$ for $n\geq 0$.
\end{Pro}

\begin{proof} Let us draw the black squares of a $2\times n$
chessboard in one $1\times n$ row, in the order we meet these squares in
the chessboard by going from
left to right. Obviously, we have a legal placement of pawns on the chessboard
if and only if we have no two
consecutive pawns in the row, or in terms of matrices and word avoidance,
the row
avoids the word \h1. The number of different legal rows is given by the
$(n+2)$-nd Fibonacci number, that is $B_{2,n}=F_{n+2}$.

\begin{center}
\begin{figure}[h]
\hspace{40pt}
\epsfxsize=300.0pt
\epsffile{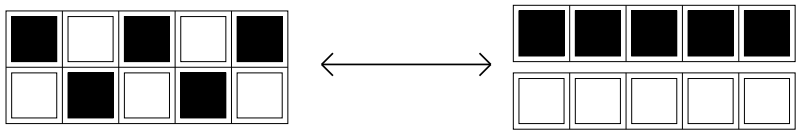}
\caption{Finding $M_{2,5}$.} \label{fig4}
\end{figure}
\end{center}

Independently, we can make the same considerations with the white squares on
the chessboard to get $W_{2,n}=F_{n+2}$. Thus, $M_{2,n}=(F_{n+2})^2$.
For instance, Figure~\ref{fig4} shows that for finding $M_{2,5}$ one can
consider two rows of length~5.
\end{proof}

\begin{Pro}\label{prop2} For all $n\geq0$, $M_{3,2n+1}=(4t_n-3t_{n-1})(2t_n-3t_{n-1})$ and $M_{3,2n}=t_n^2$, where $$t_n=\frac{1}{\sqrt{13}}\left(\left(\frac{5+\sqrt{13}}{2}\right)^{n+1}-
\left(\frac{5-\sqrt{13}}{2}\right)^{n+1}\right).$$
\end{Pro}

\begin{proof}
Let $a_n$ (resp. $b_n$) denote the number of legal placements of pawns in
the first (resp. second) shape on Figure~\ref{fig5} defined by black
squares (there are $n$ columns in each shape). According to
Observation~\ref{observ}, one has $M_{3,n}=a_nb_n$. Let us find $a_n$ and
$b_n$.

\begin{center}
\begin{figure}[h]
\hspace{60pt}
\epsfxsize=220.0pt
\epsffile{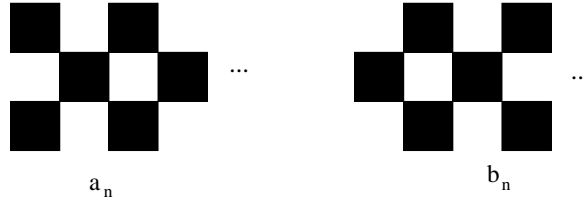}
\caption{The shapes under consideration.} \label{fig5}
\end{figure}
\end{center}

We consider the first shape. There are two black squares in the first
column. Depending on whether or not these squares have pawns, we have four
possibilities. As before, we use 1 for a square having a pawn, and 0
otherwise. Thus, the first column of the shape is either 00, or 01, or 10, or
11 when reading from top to bottom. In the first case, the first column
of the shape does not affect the rest of the shape, and therefore
can be removed. So, in the first case the number of placements of pawns is
$b_{n-1}$. In the second, third and fouth cases, the black square in the
second column must contain no pawn, that is 0, in order to have
a legal placement. This 0 does not affect what follows to the right of it,
and thus two first columns of the shape can be removed. So, the second,
third and fouth cases give $3a_{n-2}$ placements of pawns. Thus,
$a_{n}=3a_{n-2}+b_{n-1}$. Similarly, one can consider the second shape to get
$b_{n}=a_{n-1}+b_{n-2}$. Solving the equations for $a_{n}$ and $b_{n}$, we
have $$a_{2n}=t_n,\ a_{2n+1}=4t_n-3t_{n-1},\ b_{2n}=t_n,\mbox{ and }\,
b_{2n+1}=2t_n-3t_{n-1}.$$
This gives the desired result.
\end{proof}

\begin{Rem} If $a(x)$ and $b(x)$ denote the generating functions for
the numbers $a_n$ and $b_n$ respectively in the proof of Proposition~\ref{prop2}, then
$$a(x)=\frac{1+4x-3x^3}{1-5x^2+3x^4}\ \mbox{ and } \ b(x)=\frac{1+2x-3x^3}{1-5x^2+3x^4}.$$
\end{Rem}

\begin{Pro}\label{prop3} We have that $M_{4,n}={\alpha}_n^2$, where
the generating
function for the numbers ${\alpha}_n$ is given by
$$\frac{1+2x-2x^2}{1-2x-2x^2+2x^3}.$$
\end{Pro}

\begin{proof}
Let ${\alpha}_n$ (resp. ${\beta}_n$, ${\gamma}_n$, ${\delta}_n$) denote
the number of legal placements of pawns in
the first (resp. second, third, fouth) shape on Figure~\ref{fig6} defined
by black squares (there are $n$ columns in each shape). As in the proof
of Theorem~\ref{perfectSquare}, using horisontal reverse
of rows, it is easy to see that ${\alpha}_n={\beta}_n$ and
${\gamma}_n={\delta}_n$. Now, according to Observation~\ref{observ},
one has $M_{4,n}={\alpha}_n{\beta}_n=({\alpha}_n)^2$. Let us find ${\alpha}_n$.
\begin{center}
\begin{figure}[h]
\epsfxsize=380.0pt
\epsffile{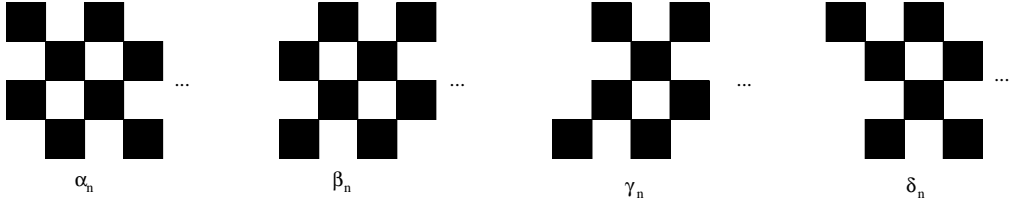}
\caption{The shapes under consideration.} \label{fig6}
\end{figure}
\end{center}
We proceed in the same way as we do in Proposition~\ref{prop2}.
If the first column of the first shape is 00 when reading the content
of the black squares from top to bottom, we can remove this column since
it does not affect the rest of the shape. So, in this case the
number of legal placements of pawns is ${\beta}_{n-1}={\alpha}_{n-1}$. If
instead of 00 we have 01 or 11, the content of the black squares in
the second column must be 00, in which case we can remove the first two
columns since they do not affect the rest of the shape. So, in this case
we have $2{\alpha}_{n-2}$ placements. The case left is when the first
column is 10. In this case the top element of the second column must be
0, and we have no information concerning the second element in this
column. Thus, in this case we have ${\gamma}_{n-1}$ replacements. Therefore,
\begin{equation}
{\alpha}_n={\alpha}_{n-1}+2{\alpha}_{n-2}+{\gamma}_{n-1}.
\label{eq_alpha}
\end{equation}

Now, to proceed further with finding ${\alpha}_n$, we need to find
${\gamma}_n$. If the element in the first column of the third shape
is 0, then we can remove this element, which gives
${\beta}_{n-1}={\alpha}_{n-1}$ replacements of pawns. If this element is
1, then the bottom element in the second column must be 0, which obviously
gives ${\delta}_{n-1}={\gamma}_{n-1}$ replacements. Thus,
\begin{equation}
{\gamma}_n={\gamma}_{n-1}+{\alpha}_{n-1}.
\label{eq_gamma}
\end{equation}
Now, from Equations~(\ref{eq_alpha}) and~(\ref{eq_gamma}) we have
$${\alpha}_n=2{\alpha}_{n-1}+2{\alpha}_{n-2}-2{\alpha}_{n-3},$$
which gives the desired result.
\end{proof}

In the way similar to that Propositions~\ref{prop2} and~\ref{prop3}
are proved, on can prove the following two propositions, which we state
without proof.

\begin{Pro}\label{prop4} For all $n\geq 0$, $M_{5,n}={\alpha}_n{\beta}_n$,
where the generating
functions for the numbers ${\alpha}_n$ and ${\beta}_n$ are given by
$$
\frac{1+7x-4x^2-7x^3+5x^4}{(1+x)(1-2x-6x^2+10x^3-4x^4)}$$
and
$$\frac{1+3x+x^2-5x^3+4x^4}{(1+x)(1-2x-6x^2+10x^3-4x^4)},$$
respectively.
\end{Pro}

\begin{Pro}\label{prop5} We have that $M_{6,n}={\alpha}_n^2$, where
the generating
function for the numbers ${\alpha}_n$ is given by
$$\frac{1+5x-9x^2-5x^3+6x^4}{1-3x-6x^2+11x^3+5x^4-6x^5}.$$
\end{Pro}


\end{document}